\documentclass[11pt]{article}
\usepackage{float} 
\usepackage[symbol]{footmisc} 
\usepackage[superscript,sort]{cite}
\usepackage{array}
\usepackage{bm}

\usepackage[normalem]{ulem} %strikeout a word
\usepackage[colorlinks=true,linkcolor=blue,citecolor=black,urlcolor=black,bookmarks=true,breaklinks=true]{hyperref}
\usepackage[usenames]{color}
\usepackage{amsmath,amssymb}
\usepackage{amsthm,amscd,amsxtra,amsfonts,amsmath,amssymb,multirow}
\usepackage{wrapfig}
\usepackage[footnotesize]{caption}
\usepackage{subcaption}
\usepackage{threeparttable} %add footnote underneath the table
\usepackage{helvet}

\usepackage{booktabs}
\usepackage[table]{xcolor}
\usepackage{xcolor,colortbl}
\usepackage{placeins}
\usepackage{graphicx}
\usepackage[T1]{fontenc}
\definecolor{MyPurple}{rgb}{1,0,1}

\providecommand{\e}[1]{\ensuremath{\times 10^{#1}}}

\newcommand{\beq}[1]{\begin{equation} \label{#1}}
\newcommand{\eeq}{\end{equation}}
\newcommand{\barray}{\begin{array}{ll}}
\newcommand{\earray}{\end{array}}

\usepackage{longtable} % longtable
\newcommand{\R}{\mathbb{R}}
\newcommand{\Z}{\mathbb{Z}}
\renewcommand{\e}{\epsilon}

\begin{document}
% \pagenumbering{roman}

% \clearpage \pagebreak \setcounter{page}{1}
% \renewcommand{\thepage}{{\arabic{page}}}

\title{Evolutionary homology on coupled dynamical systems}

\author{
Zixuan Cang$^1$, 
Elizabeth Munch$^{1,2}$
 and
Guo-Wei Wei$^{1,3,4}$ \footnote{
\{cangzixu@math, muncheli@egr, wei@math\}.msu.edu
% Zixuan Cang: cangzixu@math.msu.edu; Elizabeth Munch: muncheli@egr.msu.edu;
% Guo-Wei Wei: wei@math.msu.edu
}\\
%\address{
$^1$ Department of Mathematics \\
$^2$Department of Computational Mathematics, Science and Engineering\\
$^3$  Department of Biochemistry and Molecular Biology\\
$^4$ Department of Electrical and Computer Engineering \\
Michigan State University, MI 48824, USA \\
}

\date{\today}
\maketitle
\begin{abstract}  
% Modified by Cang
Time dependence is a universal phenomenon in nature, and a variety of mathematical models in terms of dynamical systems have been developed to understand the time-dependent behavior of real-world problems. Originally constructed to analyze the topological persistence over spatial scales,   
persistent homology has rarely been devised for time evolution. We propose the use of a new filtration function for persistent homology which takes as input the adjacent oscillator trajectories of coupled dynamical systems. We also regulate the dynamical system by a weighted graph Laplacian matrix derived from the network of interest, which embeds the topological connectivity of the network into the dynamical system. The resulting topological signatures, which we call evolutionary homology (EH) barcodes,  reveal the topology-function relationship of the network and thus give rise to the quantitative analysis of nodal properties. The proposed EH is applied to protein residue networks for protein thermal fluctuation analysis,  rendering the most accurate B-factor prediction of a set of 364 proteins. This work extends the utility of dynamical systems to the quantitative modeling and analysis of realistic physical systems.

\end{abstract}
%Key words:
%Evolutionary homology, 
%Dynamical systems 
%Protein network  

{\setcounter{tocdepth}{5} \tableofcontents}

% \newpage
\section{Introduction}

The time evolution of complex phenomena is often described by dynamical systems, i.e., mathematical models built on differential equations for continuous dynamical systems or on difference equations for discrete dynamical systems \cite{KLXia:2014b,Ott:1990, GHu:1998,Wei:2002e}. Most dynamical systems  have their origins in Newtonian mechanics. However, these mathematical models typically only admit  highly reduced descriptions of the original complex physical systems, and thus their continuous forms do not have to satisfy the Euler-Lagrange equation of the least action principle.  Although a low-dimensional dynamical system is  not expected to describe the full dynamics of a complex physical system, its long-term behavior, such as the existence of steady states (i.e., fixed points)  and/or chaotic states, offers a qualitative understanding of the underlying system.  Focused on ergodic systems, dynamic mappings, bifurcation theory, and chaos theory, the study of dynamical systems is a mathematical subject in its own right,  drawing on analysis, geometry, and topology. Dynamical systems are motivated by real-world applications, having a wide range of applications to physics, chemistry, biology, medicine, engineering, economics, and finance. Nevertheless, essentially all of the analyses in these applications are qualitative and  phenomenological in nature.

In order to pass from qualitative to quantitative evaluation of these systems, we look to the newly emerging field of topological data analysis (TDA) 
 \cite{Carlsson2009,Edelsbrunner:2010,Ghrist:2008,kaczynski:mischaikow:mrozek:04,Ghrist2014,Munch2017}. 
Homology, a tool from traditional algebraic topology, provides an algebraic structure which encodes topological structures of different dimensions in a given space, such as connected components, closed loops, and other higher dimensional analogues \cite{Hatcher}. 
To study topological invariants in a discrete data set, we use the structure of the data set, such as pairwise distance information, to build a simplicial complex, which can be loosely thought of as a generalization of a graph, and then compute the homology of the complex.
However, conventional homology is blind to translation, rotation, and scale, and thus retains too little geometric or physical information to be practically useful. 
Persistent homology, a new branch of algebraic topology,   embeds multiscale information into topological invariants to achieve an  interplay between geometry and topology.  

Given a continuum of topological spaces, called a filtration, persistence encodes the changing homology as a proxy for the shape and size of the data set by keeping track of when homological features appear and disappear over the course of the filtration.
This flexibility means that the choice of filtration allows the use of persistent homology to be tailored to the data set given and the question asked.  
As a result, it has been utilized for analysis of data sets arising from many different domains.
Just to name a few examples of particular interest to biologists\cite{Nanda2014}, 
persistence has been used in bioinformatics \cite{Kasson:2007,Perea:2015b,ZXCang:2017a,ZXCang:2017d},
neuroscience \cite{Singh:2008,Curto2008,Dabaghian:2012,curto2017can}, 
and protein folding \cite{KLXia:2014c, KLXia:2015a, KLXia:2015c,Gameiro:2014}.

The $0$-dimensional version of persistent homology was originally introduced for time series analysis and computer vision applications under the name ``size function" \cite{Fro90, Frosini:1999,Robins:1999,Robins1998,Robins2000a}.  
The generalized persistent homology theory and a practical algorithm was formulated  by Edelsbrunner {\it et al.}\cite{Edelsbrunner:2002}; the  algebraic foundation  was subsequently established by Zomorodian and Carlsson \cite{Zomorodian:2005}. 
Recently, there have been significant  developments and generalizations of the persistent homology methodology, 
\cite{BH11,CEH09,CEHM09,CGOS11,Carlsson:2009theory,CSM09,SMV11,zigzag,KLXia:2015d,Bubenik2014,OS13,DiFabio:2011}
further understanding of metrics and stability
\cite{CEH07,CCG09,cohen2010lipschitz,Bubenik2014a,Silva2017} 
and computational algorithms 
\cite{Otter2017, DFW14,Mischaikow:2013,javaPlex,Perseus, Dipha, Ripser}. 
Persistent homology is often visualized by barcodes \cite{CZOG05,Ghrist:2008} where horizontal line segments called bars represent homology generators that survive over different filtration scales. 
The equivalent persistence diagram \cite{Edelsbrunner:2010} is a representation which plots the births and deaths of the generators in a 2D plane. 
 
Persistent homology is a versatile tool for data analysis.  
However, the difficulties inherent in the interpretation of the topological space of persistence barcodes \cite{Mileyko2011,Turner2014,Munch2014} means that the most success in combining these topological signatures with machine learning methods has been found by  turning persistence barcodes into features in a well-behaved  space suitable for machine learning. 
Options for this procedure are quickly growing, and include 
persistence landscapes \cite{Bubenik:2015statistical}, 
algebraic constructions \cite{Adcock2016,CarlssonVerovsek2016,Kalisnik2018},
persistence images \cite{Adams2017,KLXia:2015c,ZXCang:2017d },
kernel methods \cite{Reininghaus2015,KLXia:2014c},
and tent functions \cite{Perea2018}. 
In 2015, Cang {\it et al.} constructed one of  the first topology based machine learning algorithms for protein classification involving tens of thousands of proteins and hundreds of tasks \cite{ZXCang:2015}.
This approach has been generalized for the predictions of protein-ligand binding affinity \cite{ZXCang:2017b} and mutation-induced protein stability change \cite{ZXCang:2017a}, and further combined with convolutional neural networks and multi-task learning algorithms \cite{ZXCang:2017c}.

The use of homology for the analysis of dynamical systems and time series analysis predates and intertwines with the beginnings of persistent homology  \cite{kaczynski:mischaikow:mrozek:04,Mischaikow1999,Gameiro2004,arai2014effects,Robins1998,Robins:1999,Robins2000a}.
More recently, there has been increased interest in the combination of persistent homology with time series analysis \cite{Robinson2014}.
Some common methods include
computing the persistent homology of the Takens embedding \cite{Perea:2015a,Perea:2015b,Perea2016,Khasawneh2015,Khasawneh2014,Khasawneh2017},  
studying the sublevelset homology of movies \cite{Kramar2015,Tralie2017}, 
and working with the additional structure afforded by persistent cohomology \cite{SMV11,Berwald2014a,Vejdemo-Johansson2015}.
Wang and Wei have defined temporal persistent homology over the solution of a partial differential equation derived from differential geometry  \cite{BaoWang:2016a}. 
This method encodes spatial connectivity into temporal persistence in the Laplace-Beltrami flow, and offers accurate quantitative prediction of fullerene isomer stability in terms of  total curvature energy for over 500  fullerene molecules. 
Closely related to our work, Stolz {\it et al.} have recently constructed persistent homology from time-dependent functional networks associated with coupled time series \cite{stolz2017persistent}. This work uses weight rank clique filtration over a defined parameter reflecting similarities between trajectories to characterize coupled dynamical systems.

The objective of the present work is to (1) define a new simplicial complex filtration using a coupled dynamical system as input, which encodes the time evolution and synchronization of the system, and (2) use the persistent homology of this filtration to study the system itself.
The resulting persistence barcode is what we call the evolutionary homology (EH) barcode. 
We are particularly interested in the encoding and decoding of the topological connectivity of a real physical system into a dynamical system.  
To this end, we regulate the dynamical system by a  generalized graph Laplacian matrix defined on a physical system with distinct topology. 
As such, the regulation encodes the topological information into the time evolution of the dynamical system. 
We use a well-studied dynamical system, the Lorenz system, to illustrate our EH formulation.   
The Lorenz attractor is utilized to facilitate the control and synchronization of chaotic oscillators by weighted graph Laplacian matrices generated from protein C$\mathrm{_\alpha}$ networks.
We create features from the EH barcodes originating from protein networks  by using the Wasserstein and bottleneck metrics.  
The resulting outputs in various topological dimensions are directly correlated with physical properties of protein residue networks. 
Finally, to demonstrate the  quantitative analysis power of the proposed EH, we apply the present method to the prediction of protein thermal fluctuations characterized by experimental B-factors.  
We show that the present EH provides some of the most accurate B-factor predictions for a set of 364 proteins. 
Our approach not only provides a new tool for quantitatively analyzing the behavior of dynamical systems but also extends the utility of dynamical systems to the quantitative modeling and prediction of important physical/biological problems.

%-----------------------------------------------------------------
%-----------------------------------------------------------------
\section{Methods}
%-----------------------------------------------------------------
%-----------------------------------------------------------------
This section is devoted to the methods and algorithms. 
In Sec.~\ref{sec:CoupledDS}, we give a brief discussion of coupled dynamical systems and their stability control via a correlation (coupling) matrix which embeds topological connectivity of a physical system into the dynamical system.  
We define persistent homology, persistence barcodes, the associated metrics, and methods for topological learning in Sec.~\ref{sec:pershom}.
We then define evolutionary homology on coupled dynamical systems in Sec.~\ref{sec:EH}.  
Finally, the full pipeline as applied to protein flexibility analysis is outlined in Sec.~\ref{ss:proteinResidueFlexibility}.    

\subsection{Coupled dynamical systems}
\label{sec:CoupledDS}
The general control of coupled dynamical systems has been well-studied in the literature. 
\cite{Ott:1990, GHu:1998,Wei:2002e,KLXia:2014b} %Don't know why, but putting this on the line above created a weird black box in the pdf.
A brief review is given to establish notation and facilitate our topological formulation, largely following the work of Hu \textit{et al.} \cite{GHu:1998}, and Xia and Wei \cite{KLXia:2014b}.  

%---
\subsubsection{Systems configuration}
We consider the coupling of $N$ $n$-dimensional dynamical systems
\begin{equation*}
\frac{d\mathbf{u}_i}{dt} = g(\mathbf{u}_i), \, i=1,2,\cdots,N,
\end{equation*}
where $\mathbf{u}_i=\{u_{i,1},u_{i,2},\cdots,u_{i,n}\}^T$ is a column vector of size $n$.
In our setup,  each $\mathbf{u}_i$ is associated to a point 
$\mathbf{r}_i \in \R^d$ which will be used to determine influence in the coupling.

The coupling of the systems can be very general, but 
% \cite{Ott:1990, GHu:1998}, including wavelets modified correlations \cite{Wei:2002e}. \liztodo{This sentece feels like it goes better in the intro}
a specific selection is an $N \times N$ graph Laplacian matrix $A$ defined for pairwise interactions
\begin{equation*}
A_{ij} = 
\begin{cases}
I(i,j),\,i\neq j,\\
-\sum\limits_{l\neq i}A_{il},\, i = j,
\end{cases}
\end{equation*}
where $I(i,j)$ is a value describing the degree of influence on the $i$th system induced by the $j$th system. %\liztodo{Does $d^{org}$ go in here somewhere?}
We assume undirected graph edges
$I(i,j) = I(j,i)$.
Let 
$\mathbf{u}=\{\mathbf{u}_1,\mathbf{u}_2,\cdots,\mathbf{u}_N\}^T$ 
be a column vector with $\mathbf{u}_i=\{u_{i,1}, u_{i,2}, \cdots, u_{i,n}\}^T$. 
The coupled system is an $N\times n$-dimensional dynamical system modeled as
\begin{equation}\label{eq:cpeq}
\frac{d\mathbf{u}}{dt} = \mathbf{G}(\mathbf{u})+\epsilon (A\otimes \Gamma)\mathbf{u},
\end{equation}
where $\mathbf{G}(\mathbf{u})=\{g(\mathbf{u}_1),g(\mathbf{u}_2),\cdots,g(\mathbf{u}_N)\}^T$, $\epsilon$ is a parameter, and $\Gamma$ is an $n \times n$ predefined linking matrix.

For this paper, we set $g$ to be the Lorenz oscillator defined as 
\begin{equation}\label{eq:Lorenz}
g(\mathbf{u}_i) =
\begin{bmatrix}
\delta(u_{i,2}-u_{i,1}) \\
u_{i,1}(\gamma-u_{i,3})-u_{i,2} \\
u_{i,1}u_{i,2}-\beta u_{i,3}
\end{bmatrix}
\end{equation}
where $\delta$, $\gamma$, and $\beta$ are parameters determining the state of the Lorenz oscillator. This system is chosen because of its relative simplicity, rich dynamics and well-understood behavior. 

%---
\subsubsection{Stability and controllability}
Let $\mathbf{s}(t)$ satisfy $d\mathbf{s}/dt = g(\mathbf{s})$.
We say the coupled systems are in synchronous state if 
\begin{equation*}
\mathbf{u}_1(t)=\mathbf{u}_2(t)=\cdots =\mathbf{u}_N(t)=\mathbf{s}(t).
\end{equation*}
The stability can be analyzed using $\mathbf{v} = \{\mathbf{u}_1-\mathbf{s}, \mathbf{u}_2-\mathbf{s}, \cdots, \mathbf{u}_N-\mathbf{s}\}^T$
% \liztodo{Isn't this the $\hat{u}$ notation defined later?}
 with the following equation obtained by linearizing Eq.~(\ref{eq:cpeq})
\begin{equation}\label{eqn:syn}
\frac{d\mathbf{v}}{dt} = [I_N\otimes Dg(\mathbf{s})+\epsilon(A\otimes\Gamma)]\mathbf{v},
\end{equation}
where $I_N$ is the $N \times N$ unit matrix and $Dg(\mathbf{s})$ is the Jacobian of $g$ on $\mathbf{s}$.% \cite{GHu:1998}. 

The stability of the synchronous state in Eq.~(\ref{eqn:syn}) can be studied by eigenvalue analysis of graph Laplacian $A$.  
Since the graph Laplacian $A$ for undirected graph is symmetric, it only admits real eigenvalues.   After diagonalizing $A$ as 
\begin{equation*}
A\phi_j = \lambda_j\phi_j, \, j=1,2,\cdots,N,
\end{equation*}
where $\lambda_j$ is the $j$th eigenvalue and $\phi_j$ is the $j$th eigenvector, $\mathbf{v}$ can be represented by 
$$\mathbf{v} = \sum\limits_{j=1}^{N}\mathbf{w}_j(t)\phi_j.$$ %\cite{GHu:1998}.
Then, the original problem on the coupled systems of dimension $N\times n$ can be studied independently on the $n$-dimensional systems %\cite{GHu:1998}
\begin{equation}\label{eq:diageq}
\frac{d\mathbf{w}_j}{dt} = (Dg(\mathbf{s})+\epsilon\lambda_j\Gamma)\mathbf{w}_j, \, j=1,2,\cdots,N.
\end{equation}
Let $L_{max}$ be the largest Lyapunov characteristic exponent of the $j$th system governed by Eq.~(\ref{eq:diageq}). 
It can be decomposed as $L_{max} = L_g + L_c$, where $L_g$ is the largest Lyapunov exponent of the system $d\mathbf{s}/dt = g(\mathbf{s})$ and $L_c$ depends on $\lambda_j$ and $\Gamma$.
In many numerical experiments in this work, we set $\Gamma=I_n$, an $n\times n$ identity matrix. 
Then the stability of the coupled systems is determined by the second largest eigenvalue $\lambda_2$. The critical coupling strength $\epsilon_0$ can, therefore, be derived as $\epsilon_0=L_g/(-\lambda_2)$. % \cite{KLXia:2014b}. 
A requirement for the coupled systems to synchronize is that $\epsilon > \epsilon_0$, while $\epsilon \leq \epsilon_0$ causes instability. 

\begin{figure}[tb]
\begin{center}
\includegraphics[width=0.8\textwidth ]{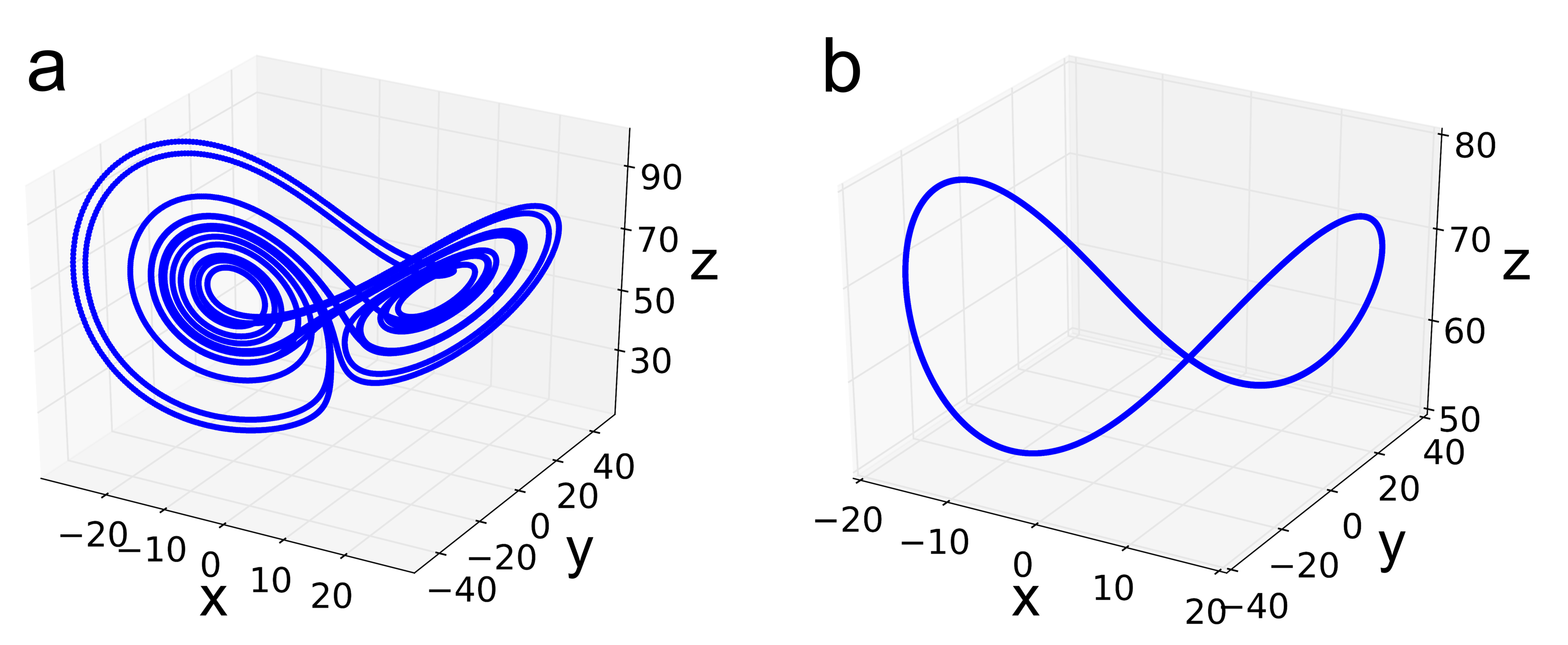}
\caption{(a) Chaotic trajectory of one oscillator without coupling. (b) The 70 synchronized oscillators associated with the carbon C$_{\alpha}$ atoms of protein PDB:1E68  are plotted together.}
\label{fig:traj}
\end{center}
\end{figure}

An example of chaos controlled by coupling is shown in Fig.~\ref{fig:traj}. 
In this example, each alpha carbon atom (C$\mathrm{_\alpha}$) of protein PDB:1E68 is associated with  a Lorenz oscillator and  the underlying locations of the oscillators are used to construct the coupling  matrix. 
The specific coupling matrix $A = A^{\mathrm{geo}} + A^{\mathrm{seq}}$ used in this example is a sum of
a graph Laplacian matrix defined using the geometric coupling,
% the oscillator geometry and a periodicity modification 
\begin{equation*}\label{eq:coupling}
A^{\rm geo}_{ij} = 
\begin{cases}
-1, \, \mathrm{if} \, i \neq j \, \mathrm{and} \, d^{\mathrm{org}}_{ij} < \epsilon_d, \\
-\sum\limits_{l\neq i}A^{\rm geo}_{il},\, i = j,
\end{cases}
\end{equation*}
and another which takes the amino acid sequence into account,
\begin{equation*}
A^{\rm seq}_{ij} = 
\begin{cases}
\epsilon_{\rm seq}, \, \mathrm{if} \, (i+1+N) \, \mathrm{mod}\, N=j, \\
-\epsilon_{\rm seq}, \, \mathrm{if} \, (i-1+N) \, \mathrm{mod}\, N=j, \\
0, \, \mathrm{otherwise}.
\end{cases}
\end{equation*}
% where $A^{\rm geo}$ cares about geometric coupling and $A^{\rm seq}$ ~ describes the coupling through an amino acid sequence. 
Here, $d^{\mathrm{org}}$ is the distance function in the original space; that is, the Euclidean distance between atoms in this example. 
% As a result, the topological connectivity of $N$ protein C$_\alpha$ atoms is embedded into $N$ Lorenz oscillators.  
The $\mathrm{mod}$ operator is used because the protein in this example is circular. 
The parameters used for the example of Fig.~\ref{fig:traj} are $\epsilon_{\rm seq}=0.7$ for sequence coupling, $\epsilon_d=4$\AA~ for spatial cutoff, and $\delta=10$, $\gamma=60$, and $\beta=8/3$ for the Lorenz system. The parameters in Eq.~(\ref{eq:cpeq}) are $\epsilon=10$ and 
\begin{equation*}
\Gamma=
\begin{bmatrix}
    0 & 0 & 0 \\
    1 & 0 & 0 \\
    0 & 0 & 0
\end{bmatrix}
.
\end{equation*}
Initial values for all  oscillators are randomly chosen.

%-----------------------
\subsection{Topological Data Analysis and Persistent Homology}
\label{sec:pershom}
%-----------------------

In this section, we introduce the necessary TDA background.  
The interested reader can find further specifics in, e.g., Carlsson \cite{Carlsson2009}, or Edelsbrunner and Harer \cite{Edelsbrunner:2010}. 

\subsubsection{Simplicial complex}

\begin{figure}[tb]
	\centering
	\includegraphics[width = 0.9\textwidth, draft = False]{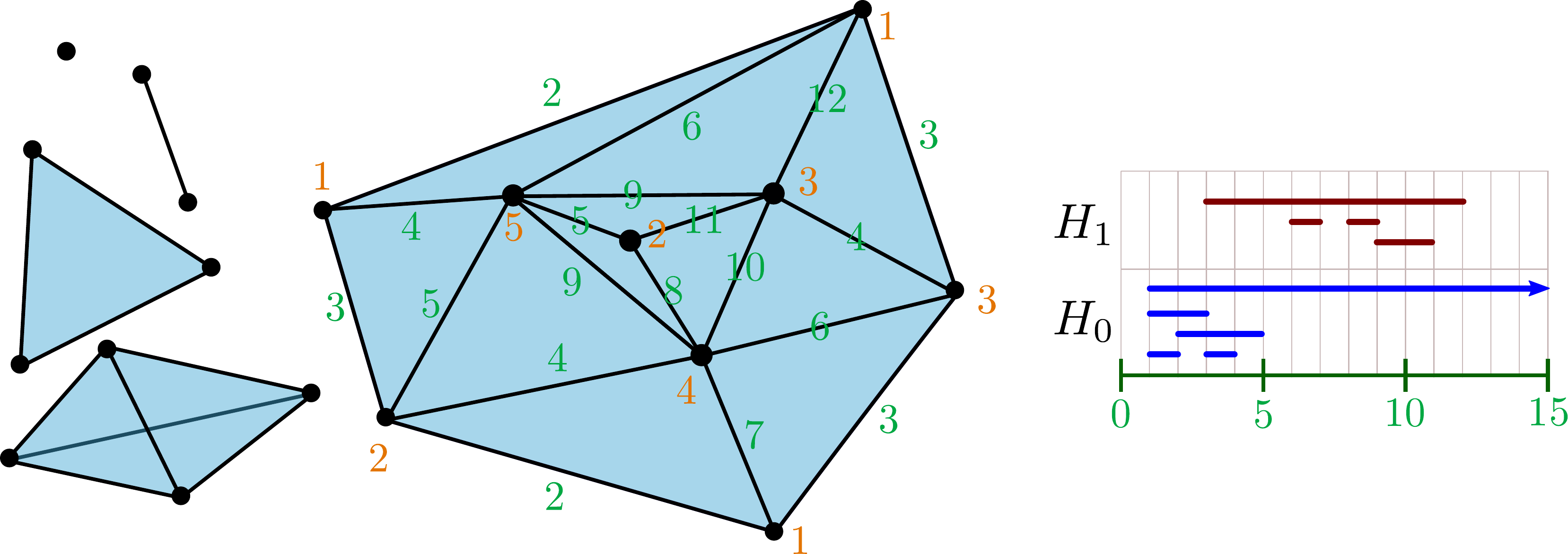}
	\caption{Examples of simplices of different dimensions (left), and a simplicial complex with a function given on the vertices and edges (middle). The barcode for the given function is drawn at right.}
	\label{fig:simplicialCpx}
\end{figure}

Topological spaces can be approximated, represented, and discretized by simplicial complexes.
An (abstract) simplicial complex is a (finite) collection of sets $K = \{\sigma_i\}_i$ where each $\sigma_i$ is a subset of a (finite) set $K^0$ called the vertex set.  
% \liztodo{Check that $N$ is used for the number of vertices}
We require that this collection satisfies the following condition: if $\sigma_i \in K$ and $\tau$ is a face of $\sigma_i$ (that is, if $\tau \subseteq \sigma_j$ commonly denoted $\tau \leq \sigma_i$), then $\tau \in K$.
If $\sigma_i$ has $k+1$ vertices, $\{v_0,v_1,\cdots,v_k\}$ where every pair of vertices is nonequivalent, $\sigma_i$ is called a $k$-simplex. 
The $k$-skeleton of a simplicial complex $K$ is the subcomplex of $K$ consisting of simplicies of dimension $k$ and below.
While the simplices for the abstract simplicial complex we will build will not have an obvious geometric meaning, there is a more geometric viewpoint from which we often reference simplicies. 
A geometric $k$-simplex can be regarded as the convex hull of $k+1$ affinely independent points in $\R^d$, and because of this we often call a $0$-simplex a point, a $1$-simplex an edge, a $2$-simplex a triangle, and a $3$-simplex a tetrahedron. 
See Fig.~\ref{fig:simplicialCpx} for an example.

\subsubsection{Homology}
The homology group for a fixed simplicial complex gives a topological characterization  which encodes holes of different dimensions. 
Homology groups are built using linear transformations called boundary operators. 
% For simplicity, we use $\mathbb{Z}_2$ coefficients. 
A $k$-chain of the simplicial complex $K$ is a finite formal sum of the $k$-simplices in $K$, $\alpha = \sum{a_i \sigma_i}$ with coefficients $a_i \in \Z_2$. 
The group of all $k$-chains with addition given by the addition of the coefficients is called the $k$-th chain group and is denoted by $C_k(K)$ or simply $C_k$ when the choice of complex is obvious. 
Note that because $\mathbb{Z}_2$ is a field, $C_k(K)$ is, in fact, a vector space.

The boundary operator $\partial_k: C_k\rightarrow C_{k-1}$ is the linear transformation generated by mapping any $k$-simplex to the sum of its codim-1 faces; namely,
\begin{equation*}
\partial_k(\{v_0,v_1,\cdots,v_k\}) = \sum\limits_{i=0}^k\{v_0,\cdots,\widehat{v}_i,\cdots,v_k\},
\end{equation*}
where $\widehat{v}_i$ means that $v_i$ is absent.
The $k$th cycle group,  $Z_k(K)$, is the kernel of the boundary operator $\partial_k$ with elements called $k$-cycles.
The $k$th boundary group, $B_k(K)$, is the image of the boundary operator $\partial_{k+1}$ and its elements are called $k$-boundaries. 
Since $\partial_k \circ\partial_{k+1}=0$, $B_k(K)$ is a subgroup of $Z_k(K)$. 
Thus we can define the $k$th homology group, $H_k(K)$,  to be the quotient group $Z_k(K)/B_k(K)$. 
% Again, note that working with field coefficients makes $H_k(K)$ a vector space.
Two $k$-cycles are called homologous if they differ by a boundary; equivalently if they are in the same equivalence class of $H_k(K)$.
Intuitively, if two $k$-cycles differ from each other by the boundary of a subcomplex, they can roughly be deformed from one to another continuously through the subcomplex. 
Each equivalence class in $H_k(K)$ can be thought of as corresponding to a $k$-dimensional ``loop'' in $K$ going around a $k+1$-dimensional ``hole'': 1-dimensional classes give information about loops going around 2D voids, 2-dimensional classes give information about enclosures of 3D voids, etc.
While the analogy is not as nice, 0-dimensional classes give information about connected components of the space.
%The dimension of $H_k(K)$ is called the $k$th Betti number.
%It is a count of $k$-dimensional holes in terms of equivalence classes.

%-------------------
\subsubsection{Filtration of a simplicial complex}

We now turn to the case where we have a changing simplicial complex and want to understand something about its structure.
Consider a finite simplicial complex $K$ and let $f$ be a real-valued function on the simplices of $K$ which satisfies the following: 
$f(\tau) \leq f(\sigma) $ for all $\tau \leq \sigma$ simplices in $K$.
We will refer to this function as the filtration function.
For any $x \in \mathbb{R}$,  the sublevelset of $K$ associated to $x$ is defined as 
\begin{equation*}\label{eq:f}
K(x) = \{\sigma\in K \mid f(\sigma) \leq x\}.
\end{equation*}
Note first that because of our assumptions on $f$,  $K(x)$ is always a simplicial complex, and second that $K(x) \subseteq K(y)$ for any $x \leq y$. 
Further, as $x$ varies, $K(x)$ only changes at the function values defined on the simplices.
Since $K$ is assumed to be finite, let $\{x_1 < x_2 < \cdots <x_\ell\}$ be the sorted range of $f$. 
The filtration of $K$ with respect to $f$ is the ordered sequence of its subcomplexes,
\begin{equation}
\label{eq:filtration}
\emptyset\subset K(x_1)\subset K(x_2) \subset \cdots \subset K(x_\ell) = K.
\end{equation}
The filtration of a simplicial complex sets the stage for a thorough topological examination of the space under multiple scales of the filtration parameter which is the output value of the filtration function $f$. 
Our choice of the filtration function $f$ for coupled dynamical systems will be given in Sec.~\ref{sec:phtime}.

\subsubsection{Persistent Homology}
The definition of homology is valid for a fixed simplicial complex, however we are interested in studying the structure of a filtration like that of Eq.~(\ref{eq:filtration}).
Functoriality of homology means that such a sequence of inclusions induces linear transformations on the sequence of vector spaces
\begin{equation}
\label{eq:Persistence}
H_k(K(x_1)) \to  H_k(K(x_2)) \to  \cdots \to  H_k(K(x_n)) .
\end{equation}
Persistent homology not only characterizes each frame in the filtration 
$\{K(x_i)\}_i$, but also tracks the appearance and disappearance (commonly referred to as births and deaths) of nontrivial homology classes as the filtration progresses.
A collection of vector spaces $\{V_i\}$ and linear transformations $f_i:V_i \to V_{i+1}$ is called a persistence module, of which Eq.~(\ref{eq:Persistence}) is an example.
It is a special case of a much more general theorem of Gabriel \cite{Gabriel1972} that sufficiently nice persistence modules can be decomposed uniquely into a finite collection of interval modules\cite{Chazal2016,Oudot2017a}.
An interval module $\mathbb{I}_{[b,d)}$ is a persistence module for which $V_i = \Z_2$ if $i \in [b,d)$ and 0 otherwise; and $f_i$ is the identity when possible, and 0 otherwise.

So, given the persistence module of Eq.~(\ref{eq:Persistence}), we can decompose it as $\bigoplus_{[b,d) \in B} \mathbb{I}_{[b,d)}$, and thus fully represent the algebraic information by the discrete collection $B$.  
These intervals exactly encode when homology classes appear and disappear in the persistence module. 
% We call this collection a persistence barcode.
The collection of such intervals can be visualized by plotting points in the 2D half plane $\{(x,y) \mid y \geq x\}$ which is known as a persistence diagram;  or by stacking the horizontal intervals, which is known as a barcode. 
In this paper, for no reason other than convenience, we represent our information using barcodes.
We call the barcode  resulting from a sequence of trivial homology groups the empty barcode and denote it by $\emptyset$.
Thus, for every interval $[b,d) \in B$, we call $b$ the birth time and $d$ the death time.

\subsubsection{Metrics on the space of barcodes}
The similarity between  persistence barcodes can be quantified by barcode space distances. 
The most commonly used metrics are the bottleneck distance \cite{CEH07} and the $p$-Wasserstein distances \cite{cohen2010lipschitz}.%,carlsson2014topological}. 
The definitions of the two distances are summarized as follows. 

The $l^{\infty}$ distance between two persistence bars $I_1 = [b_1,d_1)$ and $I_2=[b_2,d_2)$ is defined to be
\begin{equation*}
  \Delta(I_1,I_2)= \max\{ | b_2-b_1|, |d_2-d_1|\}. 
 \end{equation*}
The existence of a bar $I=[b,d)$ is measured as 
\begin{equation*}
\lambda(I):=(d-b)/2 = \min\limits_{x\in \mathbb{R}}\Delta(I, [x,x)).
\end{equation*}
This can be interpreted as measuring the smallest distance from the bar to the closest degenerate bar whose birth and death values are the same. 

% \paragraph{$\mathbf{p}$-Wasserstein distance}
For two finite barcodes 
$B_1=\{I^1_{\alpha}\}_{\alpha\in A}$ and $B_2=\{I^2_{\beta}\}_{\beta\in B}$, a  partial bijection  is defined to be a bijection $\theta:A' \to B'$ where $A'\subseteq A$ to $B'\subseteq B$.
In order to define the $p$-Wasserstein distance, we have the following penalty for $\theta$
\begin{equation*}
\label{eq:ppenalty}
P(\theta) = 
\left(
\sum\limits_{\alpha\in A'}\Delta(I^1_{\alpha},I^2_{\theta(\alpha)})^p  
+ \sum\limits_{\alpha\in A\setminus A'}\lambda(I^1_{\alpha})^p 
+ \sum\limits_{\beta \in B\setminus B'}\lambda(I^2_{\beta})^p
\right)^{1/p}.
\end{equation*}
Then the $p$-Wasserstein distance is defined as 
\begin{equation*}
d_{W,p}(B_1,B_2) = 
% \left(
\min_{\theta\in\Theta}P(\theta),
% \right)^\frac{1}{p}
\end{equation*}
where $\Theta$ is the set of all possible partial bijections from $A$ to $B$.
Intuitively, a partial bijection $\theta$ is mostly penalized for connecting two bars with large difference measured by $\Delta(\cdot)$,  and for connecting long bars to degenerate bars, measured by $\lambda(\cdot)$. 

% \paragraph{Bottleneck distance } 
The bottleneck distance is an $L_\infty$ analogue to the $p$-Wasserstein distance.
The bottleneck penalty of a partial matching $\theta$ is defined as
\begin{equation*}
\label{eq:infpenalty}
P(\theta) = \max\left\{
\max_{\alpha\in A'}\left\{
	\Delta\left(I^1_{\alpha}, I^2_{\theta(\alpha)}\right)
	\right\}, 
\max_{\alpha\in A\backslash A'}\left\{\lambda(I^1_{\alpha})\right\}, 
\max_{\beta\in B\backslash B'}\left\{\lambda(I^2_{\beta})\right\}
\right\}.
\end{equation*}
The bottleneck distance is defined as 
\begin{equation*}
d_{W,\infty}(B_1,B_2) = \min\limits_{\theta\in\Theta}P(\theta).
\end{equation*}

\subsubsection{Topological learning}
\label{sec:toplearn}

The proposed method provides a vast but relatively abstract characterization of the objects of interest. 
It is potentially powerful in quantitative analysis, but is difficult to use out of the box machine learning or data analysis techniques. 
In regression analysis or the training part of supervised learning, with 
$\mathbf{B}_i$ being the collection of sets of barcodes corresponding to the $i$th entry in the training data, the problem can be cast into the following minimization problem,
\begin{equation*}
\min\limits_{\mathbf{\theta}_{b}\in \mathbf{\Theta}_b , \mathbf{\theta}_{m}\in \mathbf{\Theta}_m} \sum\limits_{i\in I}L(\mathbf{y}_i, \mathbf{F}(\mathbf{B}_i; \mathbf{\theta}_b); \mathbf{\theta}_m),
\end{equation*}
where $L$ is a scalar loss function, %to minimize, 
$\mathbf{y}_i$ is the collection of target values in the training set, 
$\mathbf{F}$ is a function that maps barcodes to suitable input for the learning models, 
and $\mathbf{\theta}_b$ and $\mathbf{\theta}_m$ are the parameters to be optimized within the search domains $\mathbf{\Theta}_b$ and $\mathbf{\Theta}_m$ respectively. The form of the loss function also depends on the choice of metric and machine learning/regression model.

A function $\mathbf{F}$ which translates barcodes to structured representation (tensors with fixed dimension) can be used with popular machine learning models including random forest, gradient boosting trees and deep neural networks.
Another popular class of models are the kernel based models that depend on an abstract measurement of the similarity or distance between the entries. 

Our choices for $\mathbf{F}$, defined in  Eq.~(\ref{eq:phparam}) of Sec.~\ref{ss:proteinResidueFlexibility}, will arise from looking at the distance from the specified barcode to the empty barcode and there is no tuning of $\theta_b$. In Sec.~\ref{sec:bfactor} where we quantitatively analyze protein residue flexibility, we evaluate our method by checking the correlation between each topological feature defined in Eq.~(\ref{eq:phparam}) and the experimental value (blind prediction) as well as the correlation between the output of a linear regression with multiple topological features and the experimental value (regression). In the former case, there is no parameter to be optimized, while in the latter case, the specific minimization problem can be written as
\begin{equation*}
\min\limits_{\theta_m\in\mathbb{R}^{n+1}}\sum\limits_{i\in I}\left( y_i - \left[ \mathrm{EH}^{p_1,1}_i,\cdots,\mathrm{EH}^{p_n,n}_i, 1 \right]\cdot\theta_m \right)^2,
\end{equation*}
where $\mathrm{EH}^{p_k,k}_i$ is the topological parameter by computing the $p_k$-Wasserstein distance of the empty set to the $k$th barcode associated with the EH computation of the $i$th protein residue (node). $I$ is the set of indexes of all residues in the protein and $y_i$ is the experimental B-factor for the $i$th protein residue which quantitatively reflects flexibility.

%----------------------------------

\subsection{Evolutionary homology and the EH barcodes}
\label{sec:EH}

Consider a system of $N$ not yet synchronized oscillators 
$\{\mathbf{u}_1,\cdots,\mathbf{u}_N\}$ 
associated to 
a collection of $N$ embedded points, $\{\mathbf{r}_1,\cdots,\mathbf{r}_N\}\subset \R^d$.
We assume the global synchronized state is a periodic orbit denoted 
$\mathbf{s}(t)$ for $t\in[t_0,t_1]$ where 
$\mathbf{s}(t_0)=\mathbf{s}(t_1)$.
For flexibility and generality, we work on post-processed trajectories obtained by applying a transformation function on the original trajectories, $\mathbf{\widehat{u}}_i(t) := T(\mathbf{u}_i(t))$. 
The choice of function $T$ is flexible and should fit the applications; in this work, we choose 
\begin{equation}
\label{eq:trsf}
T(\mathbf{u}_i(t))=\min\limits_{t'\in[t_0,t_1]}\| \mathbf{u}_i(t)-\mathbf{s}(t')\|_2,
\end{equation}
which gives $1$-dimensional trajectories for simplicity.
Further,  in our specific example, $\mathbf{\widehat{s}}(t):=T(\mathbf{s}(t)) = 0$, but, again, this is not necessary in general.

We wish to study the effects on the synchronized system of $N$ oscillators (an $(N\times 3)$-dimensional system) after perturbing one oscillator of interest. To this end, we set the initial values of all the oscillators except that of the $i$th oscillator to $\mathbf{s}(\bar{t})$ for a fixed $\bar{t}\in [t_0,t_1]$. The initial value of the $i$th oscillator is set to $\rho(\mathbf{s}(\bar{t}))$ where $\rho$ is a predefined function playing the role of introducing perturbance to the system. After the system starts running, some oscillators will be dragged away from and then go back to the periodic orbit as the perturbance is propagated and relaxed through the system.
Let $\mathbf{\widehat{u}}^i_j(t)$ denote the modified trajectory of the $j$th oscillator after perturbing the $i$th oscillator at $t=0$. 
We focus on the subset of nodes that are affected by the perturbation,
\begin{equation*}
\label{eq:Vi}
V^i = \left\{n_j\mid \max\limits_{t>0}\{\min\limits_{t'\in[t_0,t_1]}\|\mathbf{\widehat{u}}_j^i(t)-\mathbf{\widehat{s}}(t')\|_2\} \geq \epsilon_p\right\}
\end{equation*}
for some fixed $\epsilon_p$ determining how much deviation from synchronization constitutes ``being affected''.

%-------------
\subsubsection{Filtration function defined for coupled dynamical systems}
\label{sec:phtime}
Assuming we have perturbed the oscillator for node $n_i$, let $M = |V_i|$.
We will now construct a function $f$ on the complete simplicial complex, denoted by $K$ or $K_M$, with $M$ vertices.
Here, we abuse notation and write $V_i = \{n_1,\cdots,n_M\}$.
The filtration function $f:K_M \to \R$ is built to take into account the temporal pattern of the propagation of the perturbance through the coupled systems and the relaxation (going back to synchronization) of the coupled systems. 
It requires the advance choice of three parameters: 
\begin{itemize}
	\item $\epsilon_p \geq 0$, mentioned above, which determines when a trajectory is far enough from the global synchronized state, $\mathbf{s}(t)$ to be considered unsynchronized,
	\item $\epsilon_{\mathrm{sync}}\geq 0$ which controls when two trajectories are close enough to be considered synchronized with each other, and
	\item $\epsilon_d \geq 0$ which is a distance parameter in the space where the points $\mathbf{r}_i$ are embedded, giving control on when the objects represented by the oscillators are far enough apart to be considered insignificant to each other.
\end{itemize}

We will define the function $f$ by giving its value on simplices in the order of increasing dimension.
Define
\begin{equation*}
\label{eq:tisync}
t_{\mathrm{sync}}^i = 
\min\left\{t
\mid 
\int_t^\infty \| \mathbf{\widehat{u}}_j^i(t')-\mathbf{\widehat{u}}_k^i(t')\|_2\, dt' \leq \frac{\epsilon_{\mathrm{sync}}}{2}, \, \forall j,k\right\}.
\end{equation*}
That is, $t_{\mathrm{sync}}^i$ is the first time at which all oscillators have returned to the global synchronized state after perturbing the $i$th oscillator.
The value of the filtration function for the vertex $n_j$ is defined as
\begin{equation}
\label{eq:vertexfiltration}
f(n_j) = \min
\left\{
% \min
\{ t\mid \min_{t'\in[t_0,t_1]} \|\mathbf{\widehat{u}}^i_j(t)-\mathbf{\widehat{s}}(t')\|_2 \geq \epsilon_p
\}
\cup \{t^i_{\mathrm{sync}}\}
\right\}.
\end{equation}

Next, we give the function value $f$ for the edges of $K$.
To avoid the involvement of any insignificant interaction between oscillators, an edge between $n_j$ and $n_k$ denoted by $e_{jk}$ is allowed in the earlier stage of the filtration only if 
$d^{\mathrm{org}}_{jk} \leq \epsilon_d$ 
where $d^{\mathrm{org}}_{jk}$ is the distance between $\mathbf{r}_i$ and $\mathbf{r}_j$ in $\R^d$.
Specifically, the value of the filtration function for the edge $e_{jk}$ is defined as
\begin{equation}
\label{eq:edgefiltration}
f(e_{jk}) = 
\begin{cases}
\max\left \{\min
\{t\vert \int_t^\infty \|\mathbf{\widehat{u}}^i_j(t')-\mathbf{\widehat{u}}^i_k(t')\|_2\, dt' \leq \epsilon_{\mathrm{sync}}\}, f(n_j), f(n_k) \right\}, \,
& \textrm{if} \, d_{jk}^{\mathrm{org}} \leq \epsilon_d \\
t_{\textrm{sync}}^i, \, 
& \textrm{if} \,d_{jk}^{\mathrm{org}} > \epsilon_d.
\end{cases}
\end{equation}
% where $\epsilon_{\mathrm{sync}}$ is the threshold for two trajectories to be considered as synchronized. 
It should be noted that to this point, $f$ defines a filtration function because when $d_{jk}^{\mathrm{org}}\leq\epsilon_d$, $f(n_j) \leq f(e_{jk})$ according to the definition given in Eq.~(\ref{eq:edgefiltration}). The property also holds when $d_{jk}^{\mathrm{org}} > \epsilon_d$ because $f(n_j)\leq t_{\mathrm{sync}}$ according to the definition in Eq.~(\ref{eq:vertexfiltration}) and $f(e_{jk})$ equals $t_{\mathrm{sync}}$ in this case.

\begin{figure}[tb]
\begin{center}
\includegraphics[width=0.8\textwidth, draft = False]{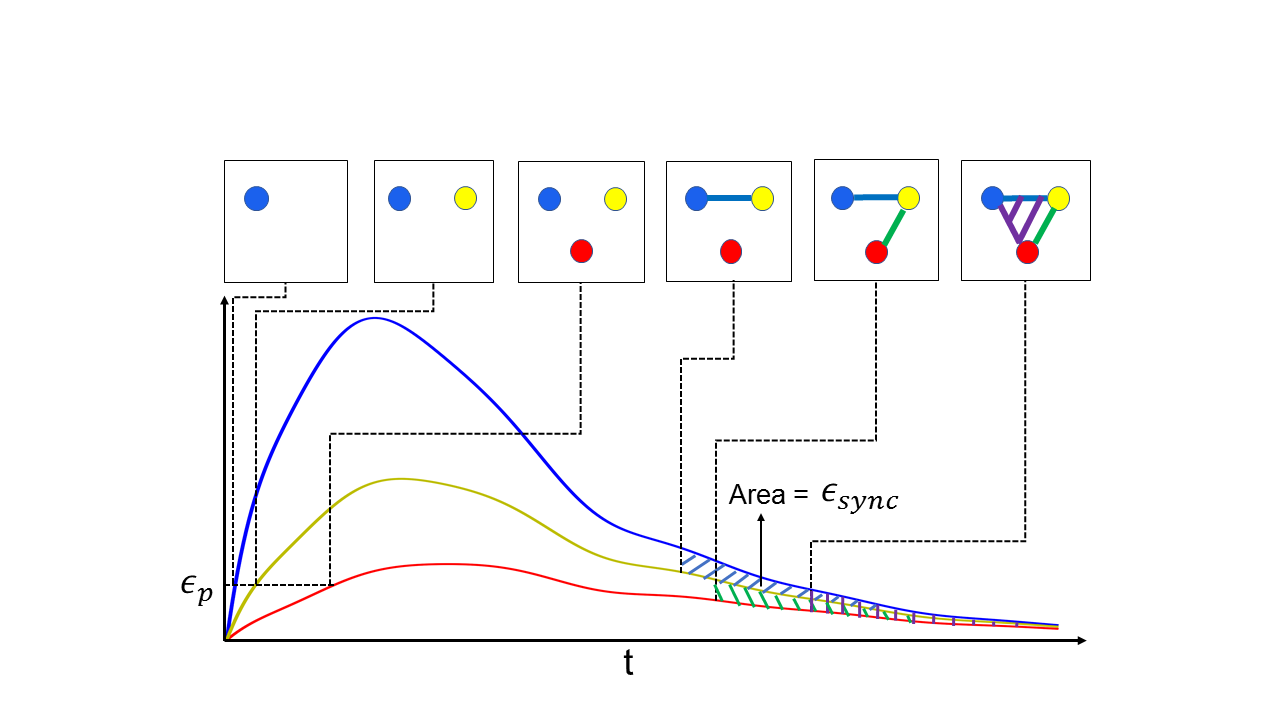}
\caption{
The filtration of the simplicial complex associated to three 1-dimensional trajectories as defined in Sec.~\ref{sec:phtime}. 
Here, each vertex corresponds to the trajectory with the same color. 
A vertex is added when its trajectory value exceeds the parameter $\e_p$; an edge is added when its two associated trajectories become close enough together that the area between the curves after that time is below the parameter $\e_{\mathrm{sync}}$.
Triangles and higher dimensional simplices are added when all necessary edges have been included in the filtration.
} 
\label{fig:filtration}
\end{center}
\end{figure}

We extend the function to the higher dimensional simplices using the definition on the 1-skeleton. 
A simplex $\sigma$ of dimension higher than one is included in $K(x)$ if all of its 1-dimensional faces are already included; that is, its filtration value is defined iteratively by dimension as
\begin{equation*}
f(\sigma) = \max_{\tau \leq \sigma} f(\tau),
\end{equation*}
where the max is taken over all codim-1 faces of $\sigma$. 
Taking the filtration of $K$ using this function (c.f.~Eq.~(\ref{eq:filtration})) means that topological changes only occur at the collection of function values $\{f(n_i)\}_i \cup \{f(e_{jk})\}_{j\neq k}$.
Fig.~\ref{fig:filtration} shows the filtration constructed for an example consisting of three trajectories.

\subsubsection{Definition of evolutionary homology}
The previous section gives a function 
% \begin{equation*}
$f_i:K_{|V^i|} \to \mathbb{R}$
% \end{equation*}
defined on the complete simplicial complex with $|V^i|$ vertices for each $i = 1,\cdots,N$.
From the filtration defined by $f_i$, we then compute the persistence barcode for homology dimension $k$, which we call the \textit{$k$th EH barcode}, denoted $B_i^k$.
The persistent homology computation for dimension $\geq1$ on the filtered simplicial complex is done using the software package Ripser \cite{Ripser} using the fact that $k$-dimensional homology only requires knowledge of $k$ and $k+1$-dimensional simplices. 
The $0$-dimensional homology is computed with a modification of the union-find algorithm.

\begin{figure}[tb]
\begin{center}
\includegraphics[width=0.75\textwidth,draft = False]{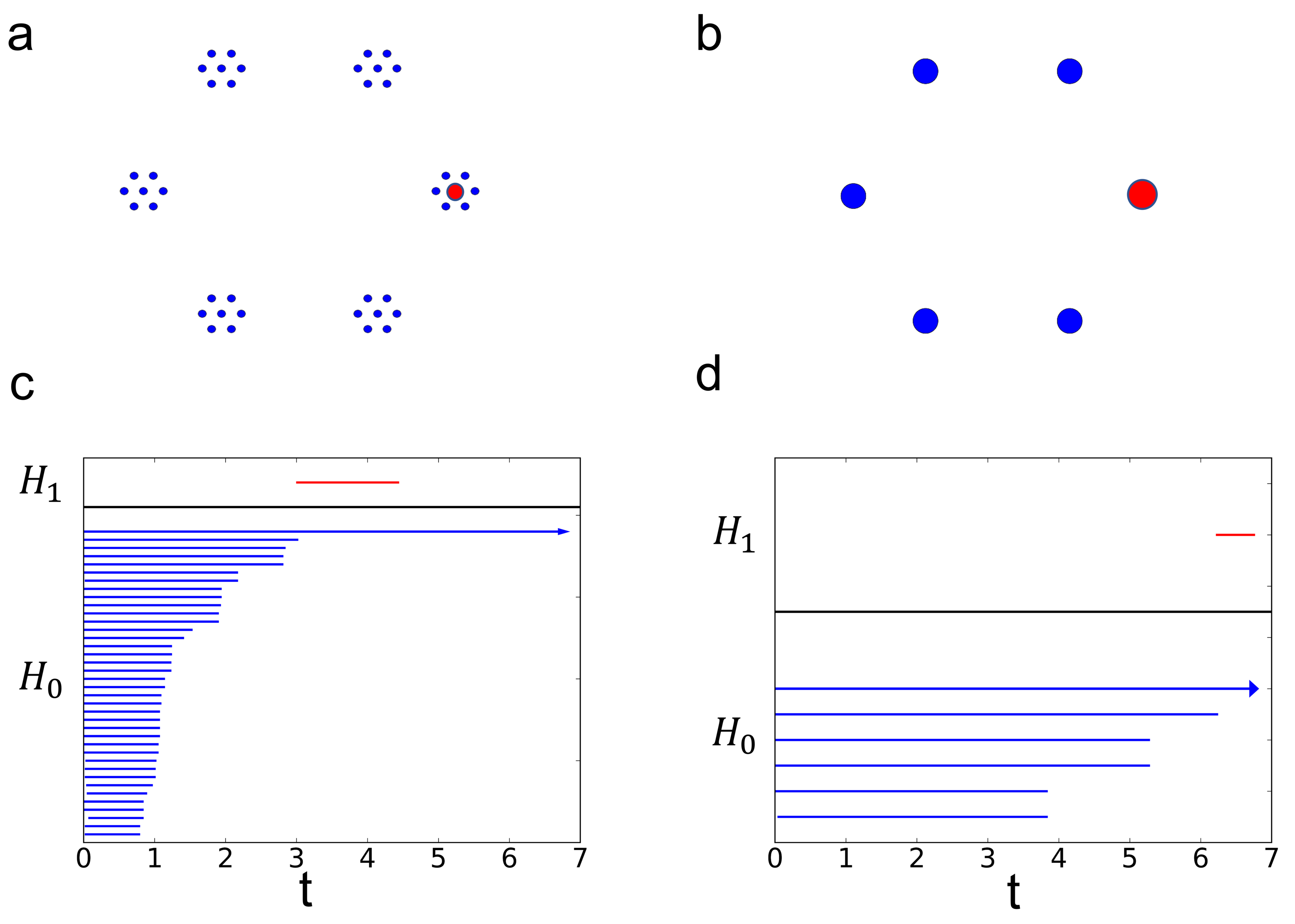}
\caption{
An example of the construction of the EH barcode. 
The geometry of two embedded systems is shown in Fig (a) and (b).
Specifically, (b) consists of six vertices of a regular hexagon with side length of $e_1$; and (a) consists of the vertices in (b) with the addition of the vertices of hexagons with a side length of $e_2\ll e_1$ centered at each of the previous vertices; 
here, $e_1 = 8$ and $e_2 = 1$.  
Figs.~(c) and (d) are the EH barcodes corresponding to Figs.~(a) and (b) respectively. 
A collection of coupled Lorenz systems is used with parameters 
$\delta=1$, 
$\gamma=12$,
$\beta=8/3$,
$\mu=8$, 
$k=2$,
$\Gamma=I_3$,
 and $\epsilon=0.12$;
see Eqs.~\eqref{eq:Lorenz}, \eqref{eq:cpmat} and \eqref{eq:cpeq}.
In the model for the $i$th residue, marked in red, the system is perturbed from the synchronized state by setting $u_{i,3}=2s_3$ with $s_3$ being the value of the third variable of the dynamical system at the synchronized state and is simulated with step size $h=0.01$ from $t=0$ using the fourth-order Runge-Kutta method. 
The calculation of persistent homology using the Vietoris-Rips filtration with Euclidean distance on the point clouds delivers similar bars corresponding to the $1$-dimensional holes in (a) and (b) which are $\left[ e_1-e_2, 2(e_1-e_2) \right)$ and $[e_1, 2e_1)$. 
}
\label{fig:Hex}
\end{center}
\end{figure}

Fig.~\ref{fig:Hex} gives an example of the geometric configurations of two sets of points associated to Lorenz oscillators and their resulting EH barcodes.  
The EH barcodes effectively examine the local properties of significant cycles in the original space which is important when the data is intrinsically discrete instead of a discrete sampling of a continuous space.
As a result, the point clouds with different geometry but similar barcodes using traditional persistence methods\footnote{Here, traditional means the Vietoris-Rips filtration on the point cloud induced by the embedding} may be distinguished by EH barcodes.

\subsection{Protein residue flexibility analysis}
\label{ss:proteinResidueFlexibility}
In this section, we combine all the methods to formulate realistic protein residue flexibility analysis using the EH barcodes. 
Consider a protein with $N$ residues and let ${\bf r}_i$ denote the position of the alpha carbon (C$_{\alpha}$) atom of the $i$th residue. 
% Each residue is associated to a Lorenz oscillator. 
The coupled systems defined in Eq.~(\ref{eq:cpeq}) are used to study protein flexibility with each protein residue represented by a $3$-dimensional Lorenz system. 
Define the distance for the atoms in the original space as the Euclidean distance between the C$\mathrm{\alpha}$ atoms, 
$d^{\mathrm{org}}({\bf r}_i,{\bf r}_j) = \|{\bf r}_i-{\bf r}_j\|_2$. 
A weighted graph Laplacian matrix is constructed based on the distance function $d^{\mathrm{org}}$ to prescribe the coupling strength between the oscillators and is defined as %an exponential function 
\begin{equation}
\label{eq:cpmat}
A_{ij} =
\begin{cases}
e^{-(d^{\mathrm{org}}({\bf r}_i, {\bf r}_j)/\mu)^\kappa}, \, i\neq j, \\
-\sum\limits_{l\neq i} A_{il}, \, i=j, 
\end{cases}
\end{equation}
where $\mu$ and $\kappa$ are tunable parameters. 
% Obviously, many other radial base functions can be used. 

\begin{figure}[tb]
\begin{center}
\includegraphics[width=0.9\textwidth]{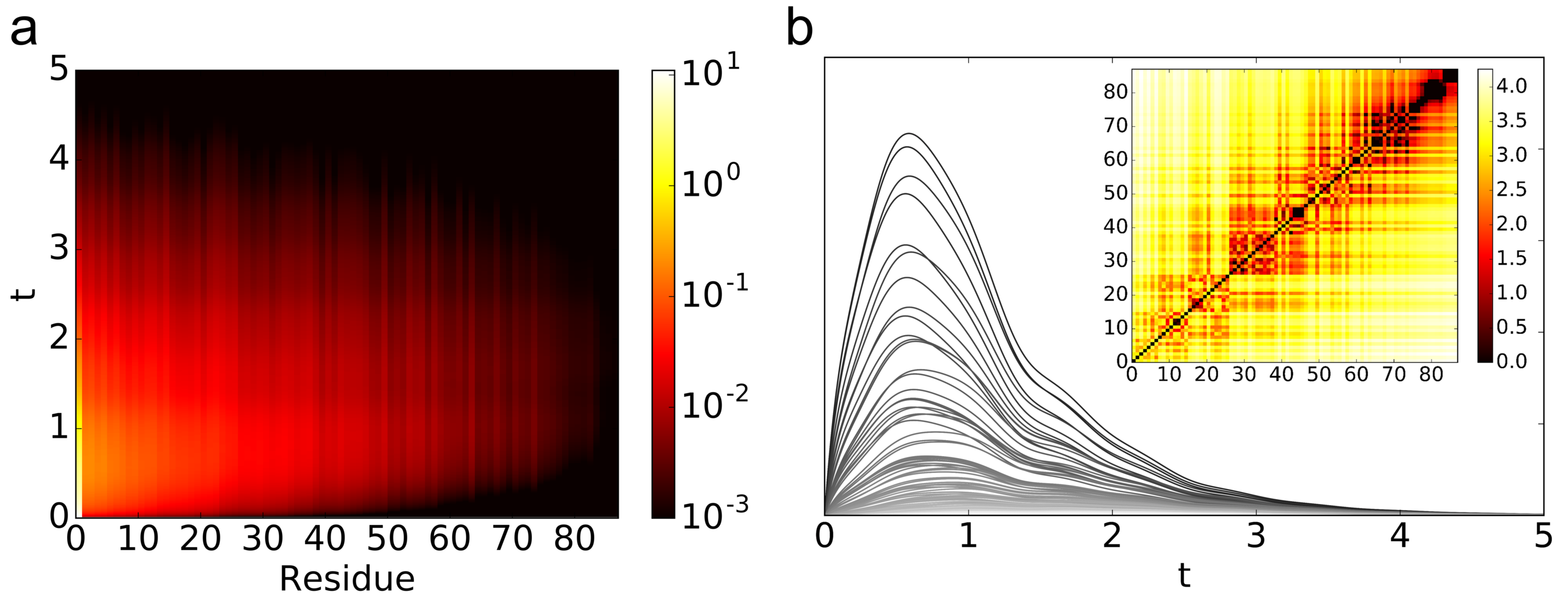}
\caption{The result of perturbing residue 31 in protein (PDB:1ABA). (a) The modified trajectories as defined in Eq.~(\ref{eq:trsf}) is plotted for each residue after the perturbation at $t=0$ as a heatmap. The residues are ordered by the (geometric) distance to the perturbed site from the closest to the farthest. (b) The modified trajectories as defined in Eq.~(\ref{eq:trsf}) is plotted for each residue after the perturbation at $t=0$ as line plots. The darker lines are closer to the perturbed site. The heatmap shows filtration value for the edges as defined in Eq.~(\ref{eq:edgefiltration}) and the order of residues is the same as in (a). 
The parameters for the coupled Lorenz system and the perturbation method are the same as that of Fig.~\ref{fig:Hex}.
}
\label{fig:heatmap}
\end{center}
\end{figure}

To quantitatively study the flexibility of a protein, one needs to extract topological information for each residue. 
To this end, we go through the process given in the previous sections once for each residue. When addressing the $i$th residue, we perturb the $i$th oscillator at a time point in a synchronized system and take this state as the initial condition for the coupled systems.
See Fig.~\ref{fig:heatmap} for an example of this procedure when perturbing the oscillator attached to a residue for a given embedding of one particular protein. 

A collection of modified trajectories $\{\mathbf{\widehat{u}}_i(t)\}_{i=1}^N$ is obtained with the transformation function defined in Eq.~(\ref{eq:trsf}). 
The persistence over time for $\{\mathbf{\widehat{u}}_i(t)\}_{i=1}^N$ is computed following the filtration procedure defined in Sec.~\ref{sec:phtime}.
Let $B_i^k$ be the $k$th EH barcode obtained from the experiment of perturbing the oscillator corresponding to residue $i$. 
We introduce the following topological features to relate to protein flexibility: 
\begin{equation}\label{eq:phparam}
\mathrm{EH}^{p,k}_i = d_{W,p}(B_i^k, \emptyset), \\
\end{equation}
where $d_{W,p}$ for $1 \leq p < \infty$ is the $p$-Wasserstein distance and $p = \infty$ is the bottleneck distance.
We will show that these features characterize the behavior of this particular collection of barcodes, which in turn, captures the topological pattern of the coupled dynamical systems arising from the underlying protein structure.

The flexibility of any given  residue is reflected by how the perturbation induced stress  is propagated and relaxed through the interaction with the neighbors. 
Such a relaxation process will induce the change in the states of the nearby oscillators. 
Therefore, the records of the time evolution of this subset of  coupled  oscillators in terms of topological invariants can be used to analyze and predict protein flexibility. 

\begin{figure}[tb]
\begin{center}
\includegraphics[width=0.9\textwidth]{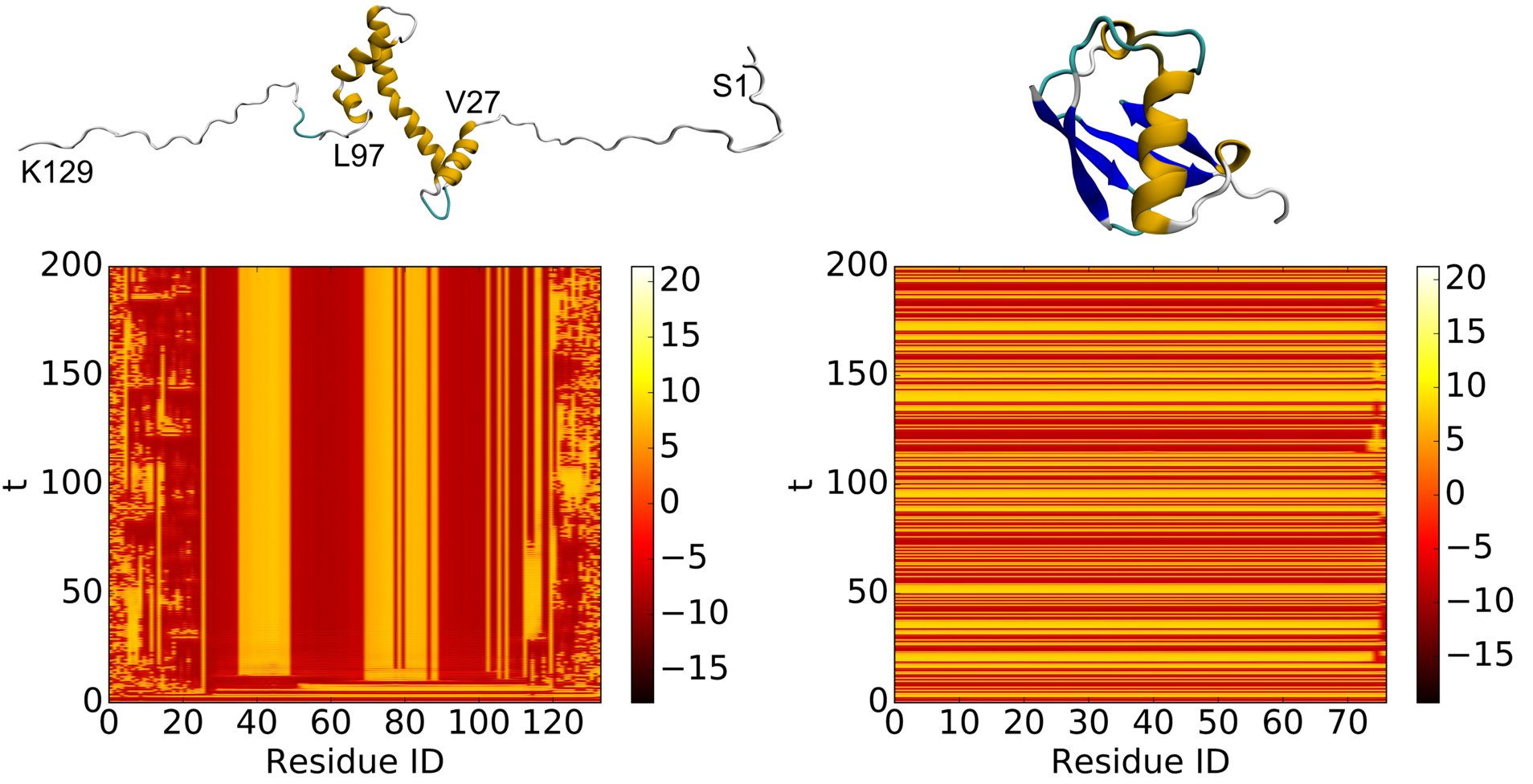}
\caption{Left: partially disordered protein, model 1 of PDB:2RVQ. Right: well folded protien, PDB:1UBQ. 
The $u_{i,1}$ value of each dynamical system is plotted as heatmap. 
The Lorenz system defined in Eq.~(\ref{eq:Lorenz}) is used with the parameters $\delta=10,\gamma=28,\beta=8/3$. 
The coupling matrix $A$ defined in Eq.~(\ref{eq:cpmat}) has parameters $\mu=14, \kappa=2$. 
The coupled system defined in Eq.~(\ref{eq:cpeq}) has parameters $\Gamma=I_3$ and $\epsilon=0.12$.
The system is initialized with a random value between $0$ and $1$ and is simulated from $t=0$ to $t=200$ with step size $h=0.01$. 
The system is numerically solved using the 4-th order Runge-Kutta method. 
It can be seen from the heatmaps that the oscillators corresponding to the disordered regions behave asynchronously.
}
\label{fig:proteinsync}
\end{center}
\end{figure}

The difference in results of the procedure can be seen in the example of Fig.~\ref{fig:proteinsync} where the control of chaotic oscillators attached to a partially disordered protein (PDB:2RVQ) and a well-folded protein (PDB:1UBQ) is demonstrated. 
Clearly, the folded part of protein 2RVQ has strong correlations or interactions among residues from residue 25 to residue 110, which leads to the synchronization of the associated chaotic oscillators. In contrast, the random coil part of   protein 2RVQ does not have much coupling or interaction among residues. Consequently,  the associated   chaotic oscillators remain in chaotic dynamics during the time evolution.   For  folded protein 1UBQ, the associated   chaotic oscillators become synchronized within a few steps of simulation, except for a small flexible tail. This behavior underpins the  use of coupled dynamical systems for protein flexibility analysis.

\section{Results}
\subsection{Disordered and flexible protein regions}\label{sec:disordered}

To illustrate the correlation between protein residue flexibility and the topological features defined in Eq.~(\ref{eq:phparam}), we study several proteins with intrinsically disordered regions. 
Intrinsically disordered proteins lack stable 3-dimensional molecular structures.
Partially disordered proteins refer to the intrinsically disordered proteins that contain both stable structure and flexible regions. In nature, the disordered regions may play important roles in biological processes which requires flexibility. 

\begin{figure}[tb]
\begin{center}
\includegraphics[width=0.9\textwidth]{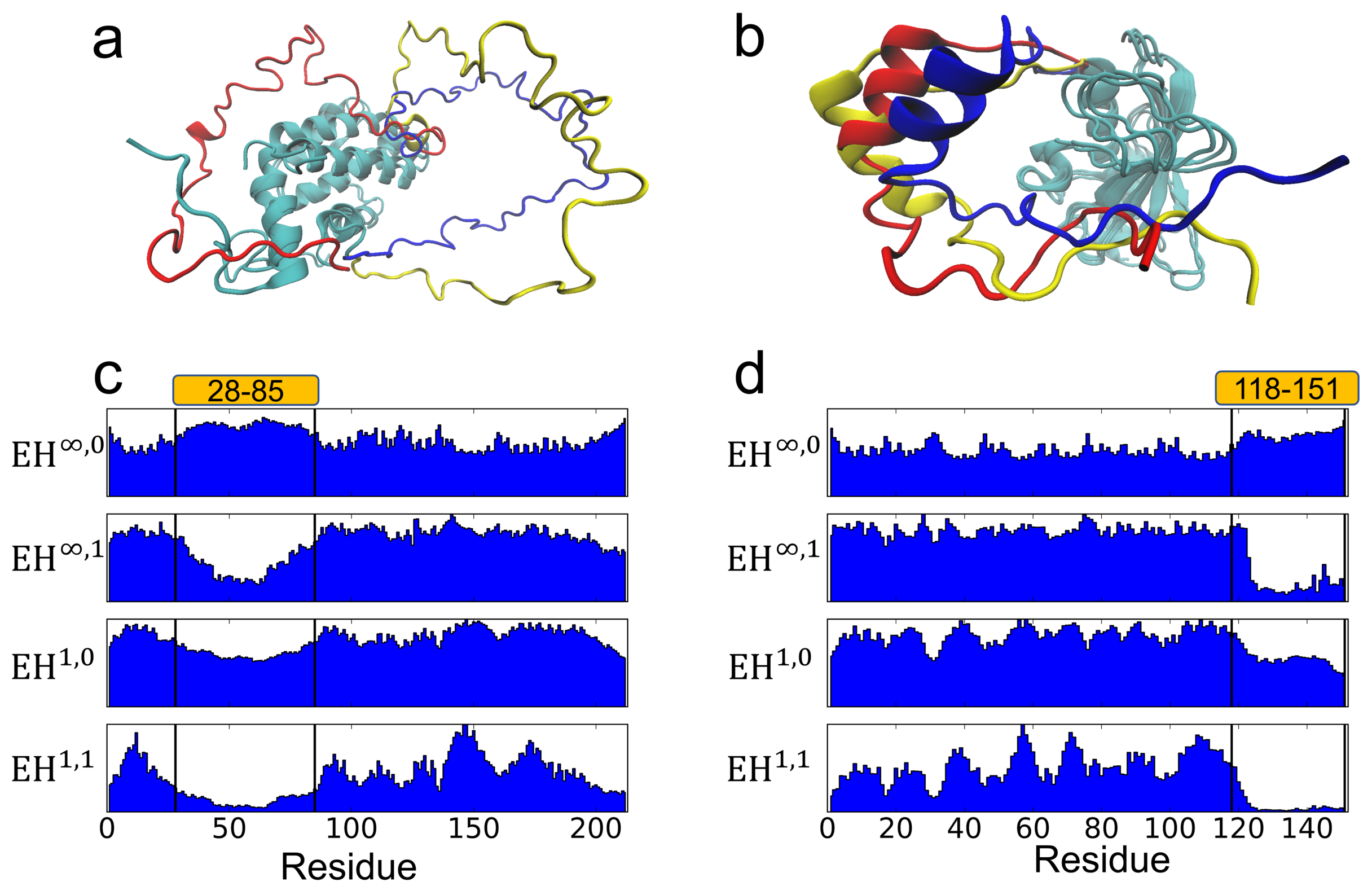}
\caption{(a) Models 1-3 of PDB:2ME9 with the disordered region colored in blue, red, and yellow for the three models. (b) Similar plot as (a) for PDB:2MT6. (c) Topological features for PDB:2ME9 whose large disordered region is from residue 28 to residue 85. (d) Topological features for PDB:2MT6 whose large disordered region is from residue 118 to residue 151. }
\label{fig:disordered}
\end{center}
\end{figure}

In what follows, we always work with the coupled Lorenz system parameters, perturbation method for the $i$th residue, and simulation described in Fig.~\ref{fig:Hex}.
The simulation is stopped when all oscillators go back to synchronized state. 
This process is repeated for each residue. 
Two NMR structures of partially disordered proteins PDB:2ME9 and PDB:2MT6 are studied. 
The topological features are computed for each model of the structures and are averaged over the models. 
The results are plotted in Fig.~\ref{fig:disordered}. 
The disordered regions clearly correlate to the peaks of EH$^{\infty,0}$ and the valleys of EH$^{\infty,1}$, EH$^{1,0}$, and EH$^{1,1}$. 
The topological features are also able to distinguish between relatively stable coils (the coils that are consistent among the NMR models) and the disordered parts (the parts that differ among the NMR models).

\subsection{Protein B-factor prediction}
\label{sec:bfactor}

Protein B-factors quantitatively measure the relative thermal motion of each atom and reflects atomic flexibility. 
The x-ray crystal structures deposited to the Protein Data Bank contain experimentally derived B-factors which can be used to validate the proposed method \cite{JKPark:2013, Opron:2014}. 
To analyze protein flexible regions, B-factor prediction is needed for protein structures built from computational models and some experimentally solved structures using NMR or cryo-EM techniques. 
Normal mode analysis (NMA) is one of the first methods proposed for B-factor predictions  \cite{Go:1983}.  
The Gaussian network model (GNM)   \cite{Bahar:1997} was known  for its better  accuracy and efficiency compared to a variety of  earlier methods \cite{LWYang:2008}. 
The multiscale flexibility-rigidity index (FRI), which is about 20\% more accurate than GNM,  has been established as the state-of-the-art in the B-factor predictions \cite{Opron:2015a}. 

\begin{figure}[tb]
\begin{center}
\includegraphics[width=0.8\textwidth]{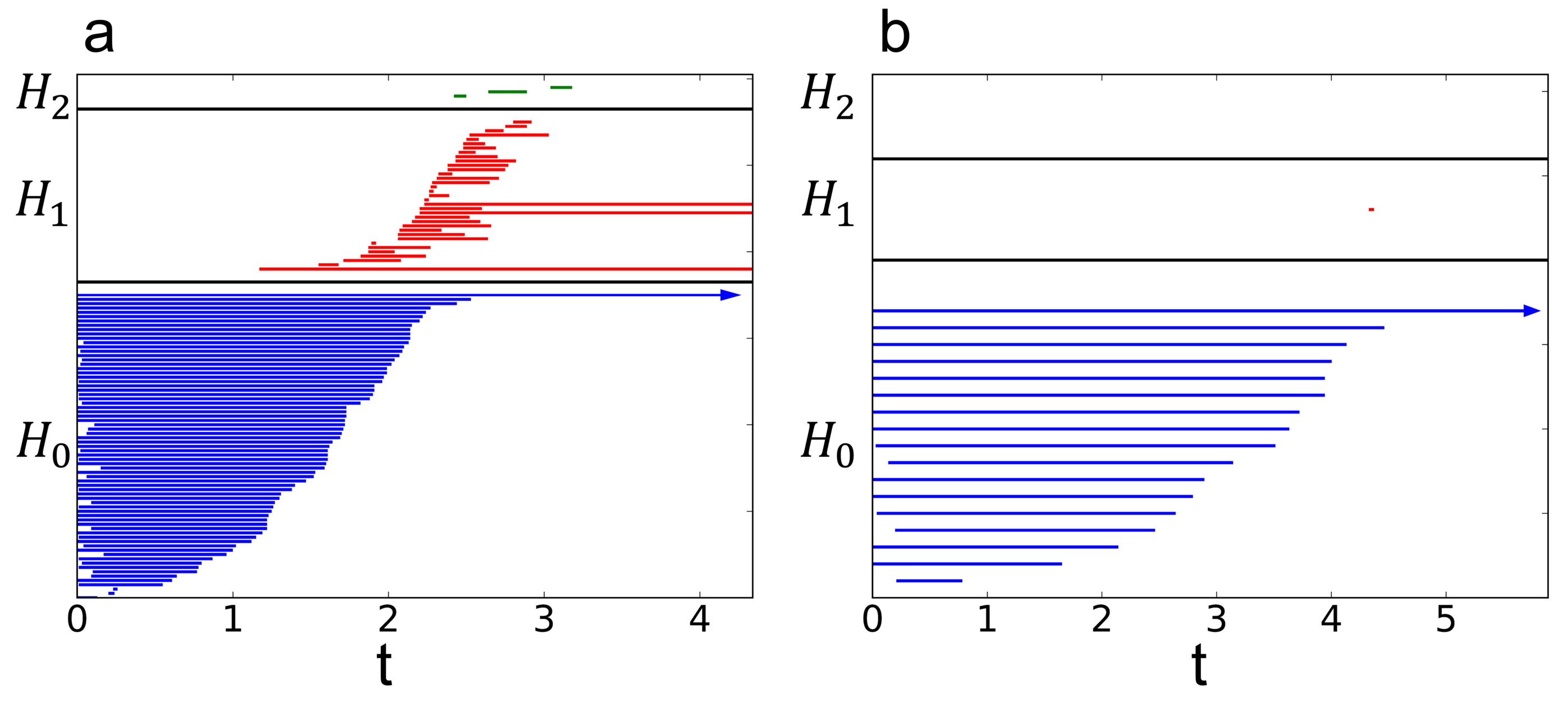}
\caption{Barcode plots for two residues. (a) Residue 6 of PDB:2NUH with a B-factor of 12.13 \AA$^2$. (b) Residue 49 of PDB:2NUH with a B-factor of 33.4 \AA$^2$.}\label{fig:barcodes}
\end{center}
\end{figure}

In this section, we compute the correlation between the topological features and the experimentally derived protein B-factors. 
We further test the proposed topological features by building a simple linear regression model with a least square penalty against the experimental B-factors. 
A collection of 364 diverse proteins reported in the literature 
% \cite{Opron:2015a} 
is chosen as the validation data (The set of 365 proteins \cite{Opron:2014} excepts PDB:1AGN due to issue in reported B-factors \cite{Opron:2015a}). 
The size of the proteins ranges from tens to thousands of amino acid residues.
The topological features in the model are the same as the setup given in Sec.~\ref{sec:disordered}. 
An example of the resulting persistence barcodes for relatively rigid and relatively flexible residues are shown in Fig.~\ref{fig:barcodes}. 
It is seen that the residue with a relatively small B-factor has many $H_0$, $H_1$ and $H_2$ bars. 
Compared to the residue having a large B-factor, it has a much richer dynamical response and barcodes with more bars.
Additionally, its $H_0$ bars are much shorter, indicating a stronger interaction with neighbor residues.

\begin{figure}[tb]
\begin{center}
\includegraphics[width=0.8\textwidth]{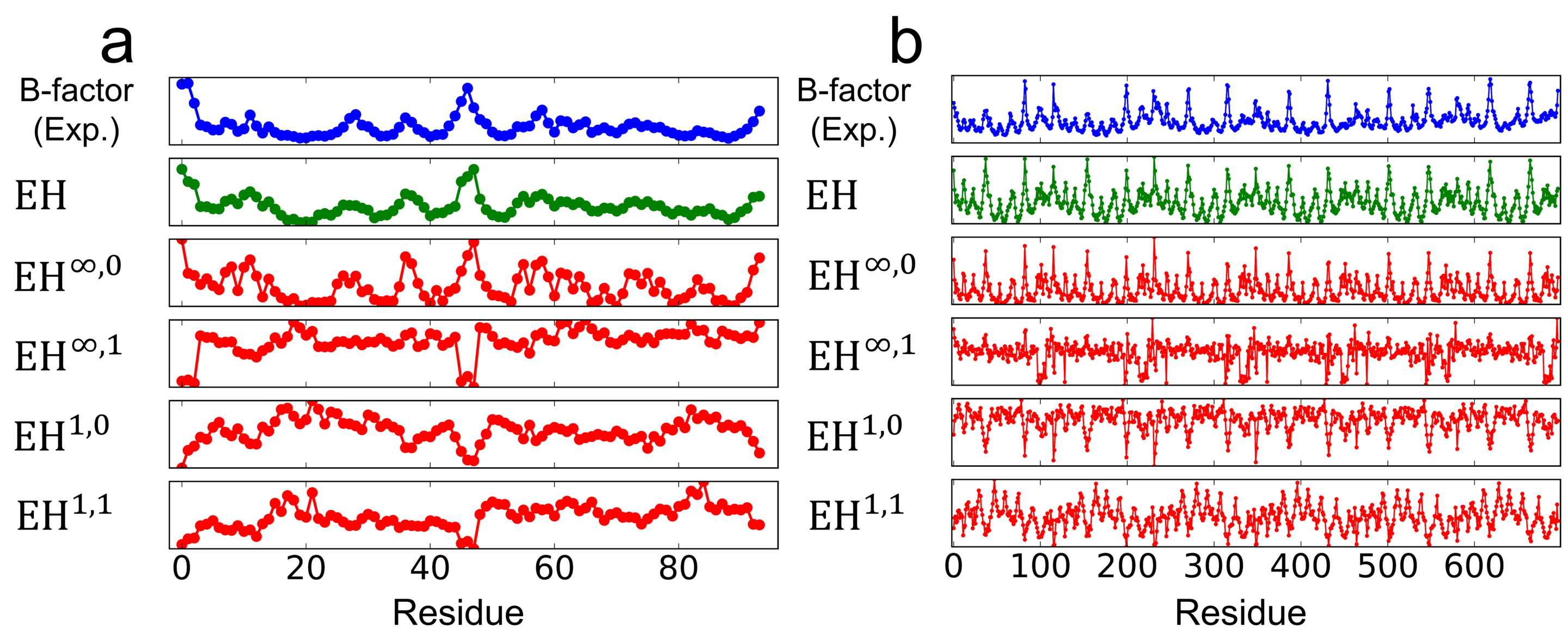}
\caption{B-factors and the computed topological features. EH shows the linear regression with EH$^{\infty,0}$, EH$^{\infty,1}$, EH$^{1,0}$, EH$^{1,1}$, EH$^{2,0}$, and EH$^{2,1}$ within each protein. (a) PDB:3PSM with 94 residues. (b) PDB:3SZH with 697 residues.}\label{fig:bfactor}
\end{center}
\end{figure}

The computed topological features are plotted against a relatively small  protein and a relatively large  protein in Fig.~\ref{fig:bfactor}. 
Clearly, 0-dimensional topological features,  specifically EH$^{\infty,0}$, provide a reasonable approximation to experimental B-factors. 
The regression using all topological information, EH, offers very good approximation to experimental B-factors.  
A summary of the results and a comparison to other methods is shown in Table \ref{tab:bfactor} for the set of 364 proteins. % \cite{Opron:2014}. 
It is seen that the present evolutionary topology based prediction outperforms  other methods in computational biophysics.  
A possible reason for this excellent performance is that the proposed method gives a more detailed description of residue interactions in terms  of three different topological dimensions and two distance metrics. 
This example indicates that the proposed EH has a great potential for other important biophysical applications, including the predictions of protein-ligand binding affinities, mutation induced protein stability changes and protein-protein interactions.

\begin{table}%\label{tab:bfactor}
\begin{center}
\rowcolors{2}{cyan!50}{white}
\begin{tabular}{|l|l|l|}
\toprule
\rowcolor{cyan!100}
Method & $R_P$ & Description \\
\midrule
EH$^{\infty,0}$ & 0.586 & Topological feature \\
EH$^{\infty,1}$ & -0.039 & Topological feature \\
EH$^{\infty,2}$ & -0.097 & Topological feature \\
EH$^{1,0}$ & -0.477 & Topological feature \\
EH$^{1,1}$ & -0.381 & Topological feature \\
EH$^{1,2}$ & -0.104 & Topological feature \\
EH$^{2,0}$ & 0.188 & Topological feature \\
EH$^{2,1}$ & -0.258& Topological feature \\
EH$^{2,2}$ & -0.100 & Topological feature \\
EH           & 0.691 & Topological features \\
% mFRI  & 0.670 & Multiscale FRI \cite{Opron:2015a} \\
% pfFRI & 0.626 & Parameter free FRI \cite{Opron:2014}  \\
% GNM  & 0.565 & Gaussian network model \cite{Opron:2014} \\
\bottomrule
\end{tabular}
\qquad
\begin{tabular}{|l|l|l|}
\toprule
\rowcolor{cyan!100}
Method & $R_P$ & Description \\
\midrule
% EH$^{B,0}$ & 0.586 & Topological metric \\
% EH$^{B,1}$ & -0.039 & Topological metric \\
% EH$^{B,2}$ & -0.097 & Topological metric \\
% EH$^{W,1,0}$ & -0.477 & Topological metric \\
% EH$^{W,1,1}$ & -0.381 & Topological metric \\
% EH$^{W,1,2}$ & -0.104 & Topological metric \\
% EH$^{W,2,0}$ & 0.188 & Topological metric \\
% EH$^{W,2,1}$ & -0.258& Topological metric \\
% EH$^{W,2,2}$ & -0.100 & Topological metric \\
EH           & 0.691 & Topological metrics \\
mFRI  & 0.670 & Multiscale FRI \cite{Opron:2015a} \\
pfFRI & 0.626 & Parameter free FRI \cite{Opron:2014}  \\
GNM  & 0.565 & Gaussian network model \cite{Opron:2014} \\
\bottomrule
\end{tabular}
\caption{The averaged Pearson correlation coefficients ($R_P$) between the computed values (blind prediction for the topological features and regression for the rest of the models) and the experimental B-factors for a set of 364 proteins \cite{Opron:2015a} (Left: Prediction $R_P$s based on EH barcodes. Right: A comparison of the $R_P$s of predictions from different  methods.). Here,  EH is the linear regression using EH$^{\infty,0}$, EH$^{\infty,1}$, EH$^{1,0}$, EH$^{1,1}$, EH$^{2,0}$, and EH$^{2,1}$ within each protein. 
For a few large and  multi-chain proteins (i.e.,  1F8R, 1H6V, 1KMM, 2D5W, 3HHP, 1QKI, and 2Q52), to reduce the computational time and as a good approximation, we compute their EH barcodes on separated (protein) chains.
We see from the table at right that the proposed EH barcode method outperforms other methods in this application.
   }
\label{tab:bfactor}
\end{center}
\end{table}

\FloatBarrier

\section{Conclusion}
Many dynamical systems are designed to understand the time-dependent phenomena in the real world. 
The topological analysis of dynamical systems is scarce in general, partially due to the fact that the topological structure of most dynamical systems is typically simple. 
In this work, we have introduced evolutionary homology (EH) to analyze the topology and its time evolution of dynamical systems. 
We present a method to embed external topology of a physical system into dynamical systems. 
EH examines the embedded topology and converts it into topological invariants over time. 
The resulting barcode representation of the topological persistence is able to unveil the quantitative  topology-function relationship of the embedded physical system. 

We have chosen the well-known  Lorenz system as an example to illustrate our EH formulation. 
An important biophysical problem, protein flexibility analysis, is employed to demonstrate the proposed topological embedding of realistic physical systems into dynamical systems. 
Specifically, we construct weighted graph Laplacian matrices from protein networks to regulate the Lorenz system, which leads to the synchronization of the chaotic oscillators associated with protein residue network nodes. 
Simplices, simplicial complexes, and homology groups are subsequently defined using the adjacent  Lorenz oscillator trajectories.  
Topological invariants and their persistence are computed  over the time evolution (filtration) of these oscillators, unveiling protein thermal fluctuations at each residue.  
The Wasserstein and bottleneck metrics are used to quantitatively discriminate  EH barcodes from different protein residues. 
The resulting model using the EH barcodes is found to outperform both geometric graph and spectral graph theory based methods in the protein B-factor predictions of a commonly used benchmark set of 364 proteins.  
      
The proposed EH method can be used to study the topological structure of a general physical system. 
Moreover, the present method extends the utility of dynamical systems, which are usually designed for qualitative analysis,  to the quantitative modeling and prediction of realistic physical systems.  
Finally, the proposed approach can be readily applied to the study of a wide variety of topology-function relationships, both within computational biology such as the role of topology in protein-ligand,  protein-protein, protein-metal and  protein-nucleic acid interactions; but also to other interactive graphs and networks in science and engineering.

\vspace{0.2in}
 
\section*{Acknowledgment} 
This work was supported in part by NSF grants  DMS-1721024  and IIS-1302285,  and MSU Center for Mathematical Molecular Biosciences Initiative.  The work of EM was supported in part by NSF grants DMS-1800446 and CMMI-1800466. 

\bibliographystyle{unsrt}
\bibliography{refs_Short}

\end{document}